\documentstyle[12pt,amssymb]{article}
\begin{document}
\newtheorem{Theo}{Theorem}
\newtheorem{Ex}{Example}
\newtheorem{Def}{Definition}
\newtheorem{Lem}{Lemma}
\newtheorem{Cor}{Corollary}
\newenvironment{pf}{{\noindent\bf Proof.\ }}{\ $\Box$\medskip}
\newtheorem{Prop}{Proposition}
\renewcommand{\theequation}{\thesection.\arabic{equation}}

%%%%%%%%%%%%%%%%%%%%%%%%%%%%%%%%%%%%%%%%%%%%%%%%%%%%%%%%%%%%%%%%%

\mathchardef\za="710B  %\alpha
\mathchardef\zb="710C  %\beta
\mathchardef\zg="710D  %\gamma
\mathchardef\zd="710E  %\delta
\mathchardef\zve="710F %\epsilon
\mathchardef\zz="7110  %\zeta
\mathchardef\zh="7111  %\eta
\mathchardef\zvy="7112 %\theta
\mathchardef\zi="7113  %\iota
\mathchardef\zk="7114  %\kappa
\mathchardef\zl="7115  %\lambda
\mathchardef\zm="7116  %\mu
\mathchardef\zn="7117  %\nu
\mathchardef\zx="7118  %\xi
\mathchardef\zp="7119  %\pi
\mathchardef\zr="711A  %\rho
\mathchardef\zs="711B  %\sigma
\mathchardef\zt="711C  %\tau
\mathchardef\zu="711D  %\upsilon
\mathchardef\zvf="711E %\phi
\mathchardef\zq="711F  %\chi
\mathchardef\zc="7120  %\psi
\mathchardef\zw="7121  %\omega
\mathchardef\ze="7122  %\varepsilon
\mathchardef\zy="7123  %\vartheta
\mathchardef\zf="7124  %\varomega
\mathchardef\zvr="7125 %\varrho
\mathchardef\zvs="7126 %\varsigma
\mathchardef\zf="7127  %\varphi
\mathchardef\zG="7000  %\Gamma
\mathchardef\zD="7001  %\Delta
\mathchardef\zY="7002  %\Theta
\mathchardef\zL="7003  %\Lambda
\mathchardef\zX="7004  %\Xi
\mathchardef\zP="7005  %\Pi
\mathchardef\zS="7006  %\Sigma
\mathchardef\zU="7007  %\Upsilon
\mathchardef\zF="7008  %\Phi
\mathchardef\zW="700A  %\Omega

\newcommand{\be}{\begin{equation}}
\newcommand{\ee}{\end{equation}}
\newcommand{\lra}{\longrightarrow}
\newcommand{\ra}{\rightarrow}
\newcommand{\bea}{\begin{eqnarray}}
\newcommand{\eea}{\end{eqnarray}}
\newcommand{\beas}{\begin{eqnarray*}}
\newcommand{\eeas}{\end{eqnarray*}}
\newcommand{\Z}{{\Bbb Z}}
\newcommand{\R}{{\Bbb R}}
\newcommand{\C}{{\Bbb C}}
\newcommand{\Li}{{\cal L}}
\newcommand{\SL}{SL(2,\C)}
\newcommand{\Sl}{sl(2,\C)}
\newcommand{\SU}{SU(2)}
\newcommand{\su}{su(2)}
\newcommand{\SB}{SB(2,\C)}
\newcommand{\Sb}{sb(2,\C)}
\newcommand{\G}{{\cal G}}
\newcommand{\g}{{\frak g}}
\newcommand{\D}{{\rm d}}
\newcommand{\de}{\,{\stackrel{\rm def}{=}}\,}
\newcommand{\we}{\wedge}
\newcommand{\nn}{\nonumber}
\newcommand{\ot}{\otimes}
\newcommand{\s}{{\textstyle *}}
\newcommand{\ts}{T^\s}
%%%%%%%%%%%%%%%%%%%%%%%%%%%%%%%%%%%%%%%%%%%%%%%%%%%%%%%%%%%%%%%%%%
\newcommand{\pa}{\partial}
\newcommand{\ti}{\times}
%%%%%%%%%%%%%%%%%%%%%%%%%%%%%%%%%%%%%%%%%%%%%%%%%%%%%%%%%%%%%%%%%%
%\begin{document}
\begin{center}
{\Large\bf Generalized $n$-Poisson brackets on a symplectic manifold}
\vskip 1cm
J. Grabowski\footnote{Institute of Mathematics, Warsaw University,
ul. Banacha 2, 02-097 Warszawa, Poland;
{\it e-mail:} jagrab@mimuw.edu.pl .\\
Supported by KBN, grant No. 2 P03A 042 10.} \\ and \\
G. Marmo\footnote{
Dipartimento di Scienze Fisiche, Universit\`a di Napoli,
Mostra d'Oltremare, Pad. 20, 80125 Napoli, Italy;
{\it e-mail:} gimarmo@na.infn.it .\\
This work has been partially supported by PRIN-97 "SINTESI".}
\end{center}
\date{\ }
\centerline{\bf Abstract}
\bigskip\noindent
On  a  symplectic  manifold  $(M,\omega)$  we study a  family  of
generalized Poisson brackets associated with  $2k$-forms  $\omega^k$.  The
extreme  cases  are  related  to  the  Hamiltonian  and  Liouville
dynamics. We show that the Dirac brackets  can  be  obtained  in  a
similar way.
\bigskip\noindent
\setcounter{equation}{0}
\section{Introduction}
Hamiltonian formalism and Poisson brackets have acquired a dominant role
in the description of classical systems after their use by Dirac in the
formulation of Quantum Mechanics \cite{Di1}.
However, when dealing with statistical mechanics on the phase space, the
Liouville measure plays a more relevant role.
On a symplectic manifold (phase space), as noticed  by  Poincar\'e
\cite{Po}, there are also available other integral invariants.
To be more specific, on any
$2n$-dimensional symplectic manifold $(M,\zw)$ any Hamiltonian system in
addition to $\zw$ preserves also $\zw^2,\zw^3,\dots,\zw^n$.
Among them $\zw$ and $\zw^n$ play a privileged role, because they define
isomorphisms between covariant and contravariant tensors, intermediate
powers do not.
\par
Vector fields preserving a volume form are  divergenceless  vector
fields, i.e. they
represent continuity equations and therefore define `conserved
quantities'. Vector fields preserving a volume form have been called {\em
Liouville dynamical systems} \cite{MSSZ}.
They furnish a geometrical approach to all
dynamical systems which satisfy some continuity conditions.
For these systems it is possible to study the analogue of Poisson brackets.
\par
In this paper we would like to show that it is possible to introduce and
study the analogue of Poisson brackets also for intermediate powers of $\zw$,
hopefully the analysis of these situations from the dynamical point of view
may bring in a finer classification of dynamical systems which goes
beyond the dichotomy Hamiltonian--non-Hamiltonian dynamics.
\par
We introduce very briefly what is the ideology in introducing brackets
associated with powers of $\zw$.
The main observation is that, if on a manifold $M$
of dimension $m$ we have a volume element ($m$-form) $\zW$ and an
$(m-2)$-form $\za$, for any two functions $f,g\in C^\infty(M)$ we can
define the bracket $\{ f,g\}$ by setting
\be
\{ f,g\}\zW=\D f\we\D g\we\za.
\ee
Of course, we have to put additional requirements on $\za$, if  we
want to have the Jacobi identity for the bracket.
On a symplectic manifold $(M,\zw)$ we recover the standard Poisson bracket
if we put $\zW=\zw^n$ and $\za=n\zw^{n-1}$.
\par
If we use an $(m-4)$-form $\zb$, we can define a quaternary bracket by
\be
\{ f_1,f_2,f_3,f_4\}\zW=\D f_1\we\D f_2\we\D f_3\we\D f_4\we\zb
\ee
and so on. When the form $\zg$ is just a function, we get
\be
\{ f_1,\dots,f_m\}\zW=\zg\D f_1\we\cdots\we\D f_m
\ee
or
\be
\{ f_1,\dots,f_m\}=\zg{\rm det}\left(\frac{\pa f_i}{\pa x_j}\right).
\ee
which is the prototype of a  Nambu  bracket  \cite{Na}  (cf.  also
\cite{GMP}).
\par
Previous idea can be used also to deal with brackets in the presence
of `second class constraints' \cite{Di1,Di2,MMS}.
If $\zq_1,\zq_2,\dots,\zq_{2k}$ are functions on a symplectic
manifold such that
\bea
&\D\zq_1\we\D\zq_2\we\cdots\we\D\zq_{2k}\ne 0,\\
&{\rm det}\left(\{\zq_j,\zq_k\}\right)\ne 0,
\eea
then we may define the `Dirac bracket' $\{\ ,\ \}_D$ by
\bea\nn
&(\D f\we\D g)\we\D\zq_1\we\zq_2\we\cdots\we\D\zq_{2k}\we\zw^{n-k-1}=\\
&\{ f,g\}_D\D\zq_1\we\zq_2\we\cdots\we\D\zq_{2k}\we\zw^{n-k}.
\eea

Various generalizations of Poisson brackets, recently dealt within the
literature, admit as a prototype those associated with powers of the
symplectic structure on a symplectic manifold or a variant of it.
These brackets are usually expressed in terms of multivector fields.
The multivector fields giving rise to our brackets are defined  by
$i_{\Lambda_\omega}\Omega=\omega$ (see our Theorem 1  in  the  next
section).  For  bivector  fields   we   recall   \cite{GMP}   that
$$i_{[\Lambda_1,\Lambda_2]}\Omega=-i_{\Lambda_1}i_{\Lambda_2}{\rm
d}\Omega-{\rm d}i_{\Lambda_2\wedge\Lambda_1}\Omega  +i_{\Lambda_1}
{\rm                       d}i_{\Lambda_2}\Omega+i_{\Lambda_2}{\rm
d}i_{\Lambda_1}\Omega.$$
A bivector field $\Lambda$ defines a Poisson bracket if  and  only
if
$${\rm    d}i_{\Lambda\wedge\Lambda}\Omega=2i_{\Lambda}{\rm     d}
i_\Lambda\Omega.$$
In this way we can deal
also with odd-dimensional manifolds (contact manifolds, for instance).
The generalization, to include also Jacobi brackets, requires the
introduction of brackets on modules rather than on rings of functions.
A manifold $M$ which is furnished with a bi-vector field $\zL$ and a vector
field $X$ satisfying
\bea
&[X,\zL]=0,\\
&[\zL,\zL]=2X\we\zL,
\eea
where the brackets are the Schouten brackets, is called
{\em Jacobi manifold} with the {\em Jacobi bracket} \cite{DLM}
\be
\{ f,g\}=\zL(f,g)+fX(g)-gX(f).
\ee
If now, on the manifold $M\ti\R$, we consider the bracket associated to
the bivector field $e^{-2s}(\zL+\pa_s\we X)$ and evaluate it on the
$C^\infty(M)$-module of functions
$\{ \tilde f=e^sf: f\in C^\infty(M)\}$, we find
that we recover the Jacobi bracket on $M$:
\be
e^{-2s}(\zL+\pa_s\we X)(\tilde f,\tilde g)=\zL(f,g)+fX(g)-gX(f).
\ee

\section{Generalized $n$-Poisson brackets}

The following  theorem  describes  the  relation  of  $k$-brackets
defined  by  differential  forms  (a  volume   $m$-form   and   an
$(m-k)$-form) to multivector fields.

\begin{Theo} Let $\zW$ be a volume $m$-form on a manifold $M$
and let $\za$ be an $(m-k)$-form. Then the $k$-bracket of functions
defined by
\be
\{ f_1,\dots,f_k\}\zW=\D f_1\we\cdots\we\D f_k\we\za
\ee
is generated by the $k$-vector field $\zL$ defined by $i_\zL\zW=\za,$
i.e.
\be
\{ f_1,\dots,f_k\}=<\D f_1\we\cdots\we\D f_k,\zL>.
\ee
\end{Theo}
\begin{pf} The bracket satisfies clearly the Leibniz rule, so it is
generated by a $k$-vector field $\zL$. Contractions with the volume
form $\zW$ give rise to isomorphisms between the corresponding contravariant
and covariant tensors and we have just to prove that
\be\label{e1}
<\D f_1\we\cdots\we\D f_k,\zL>\zW=\D f_1\we\cdots\we\D f_k\we(i_\zL\zW).
\ee
Since it is enough to prove (\ref{e1}) pointwise, we may just work in
a vector space $V$ with a basis $X_1,\dots,X_m$ and the volume
$\zW=X_1\we\cdots\we X_m$. Let $X^*_1,\dots,X^*_m$ be the dual basis.
We have to prove that for any $I=(i_1,\dots,i_k)$, $1\le i_1<\cdots<
i_k\le m$, we have
\be\label{e2}
<X_{i_1}\we\cdots\we X_{i_k},\zL>\zW=
X_{i_1}\we\cdots\we X_{i_k}\we(i_\zL\zW)
\ee
for any $\zL\in\zL^kV^*.$ Let us write
\be
\zL=\sum_{I=(i_1,\cdots,i_k)}a_IX_{i_1}^*\we\cdots\we X_{i_k}^*.
\ee
Then
\be
i_\zL\zW=\sum_{I=(i_1,\cdots,i_k)}a_Ii_{X_{i_1}^*\we\cdots\we X_{i_k}^*}
X_1\we\cdots\we X_m
\ee
and (\ref{e2}) reduces to
\be
X_{i_1}\we\cdots\we X_{i_k}\we i_{X_{i_1}^*\we\cdots\we X_{i_k}^*}
(X_1\we\cdots\we X_m),
\ee
which is obvious.
\end{pf}

\medskip\noindent
\begin{Lem} If $\zL_1,\dots,\zL_m$ are pairwise compatible Poisson tensors,
i.e. $[\zL_i,\zL_j]=0$, where $[\ ,\ ]$ is the Schouten bracket,
then any wedge-product of the tensors $\zL_1,\dots,\zL_m$ commutes with
any other wedge-product of them with respect to the Schouten bracket.
In particular, if $\zL$ is a Poisson tensor, then
\be\label{p}
[\zL^i,\zL^j]=0
\ee
for all $i,j=1,2,\dots$, where
$$\zL^i=\underbrace{\zL\we\cdots\we\zL}_
{i-{\rm times}}.$$
\end{Lem}
\begin{pf}
It follows immediately from the Leibniz rule for the Schouten bracket
\end{pf}

\medskip\noindent
{\bf Remark.} The same remains valid for any Nambu-Poisson structure
$\zL$ by similar arguments.

\medskip\noindent
Let now $\zw$ be a symplectic form on an $2n$-dimensional manifold $M$ and
let $\zL$ be the corresponding Poisson tensor $\zL=\zw^{-1}$.
\begin{Lem}
\be
i_\zL\zw^k=k(n-k+1)\zw^{k-1}.
\ee
\end{Lem}
\begin{pf} Working in a Darboux chart, we have
\beas
&i_\zL\zw^k=i_{(\sum_{j=1}^n\pa_{p_j}\we\pa{q_j})}\zw^k=
\sum_{j=1}^ni_{\pa_{q_j}}i_{\pa_{p_j}}\zw^k=\\
&\sum_{j=1}^ni_{\pa_{q_j}}(k\D q_j\we\zw^{k-1})=
k\sum_{j=1}^n(\zw^{k-1}+(k-1)\D q_j\we\D p_j\we\zw^{k-2})=\\
&kn\zw^{k-1}-k(k-1)\zw^{k-1}=k(n-k+1)\zw^{k-1}.
\eeas
\end{pf}

\begin{Theo} The $2k$-bracket defined by
\be\label{e3}
\{ f_1,\cdots,f_{2k}\}\frac{\zw^n}{n!}=
k!\D f_1\we\cdots\we\D f_{2k}\we\frac{\zw^{n-k}}{(n-k)!}
\ee
is generated by the $2k$-vector field $\zL^k$. It is an $2k$-Poisson
bracket in the sense of \cite{APP1, APP2, APP3}.
\end{Theo}
\begin{pf} Using Lemma 2 we can prove inductively that
\be
i_{\zL^k}\frac{\zw^n}{n!}=\underbrace{i_\zL\cdots i_\zL}_
{k-{\rm times}}\frac{\zw^n}{n!}=k!\frac{\zw^{n-k}}{(n-k)!}
\ee
which shows, in view of Lemma 2, that the bracket is induced by $\zL^k$.
According to Lemma 2, $\zL^k$ is an $2k$-Poisson structure.
\end{pf}

{\bf Example} We shall consider $M=T^*\R^3$  with  the  symplectic
structure  $\zw_B=\zw_0+\ze_{ijk}B^i\D  q^j\we\D   q^k$.   It   is
possible to compute $\zw^2$ and $\zw^3$ and find
\bea
\zw^2&=&\zw_0^2+(B^1\D     p_1+B^2\D      p_2+B^3\D      p_3)\we\D
q^1\we\D q^2\we\D q^3,\\
\zw^3&=&\zw_0^3.
\eea
It is interesting to notice that $\zw^3$ does not  keep  track  of
the presence of the magnetic field.  This  property  is  sometimes
quoted to account for the fact that there is no `diamagnetism'  at
the classical level.
\par
The use of $\zw^2$ for computing brackets in the form
\be
\{ f_1,f_2\}\zW=\D f_1\we\D f_2\we\zw^2
\ee
will reproduce  the  standard  bracket  associated  with  $\zw_B$,
however in this way of computing we show immediately  that,  while
$\{ q^i,q^j\}=0$, we find now $\{ p_i,p_j\}=\ze_{ijk}B^k$ and  the
Jacobi identity is equivalent to
\be
{\rm div}\vec{B}=\frac{\pa B^1}{\pa q^1}+\frac{\pa B^2}{\pa q^2}+
\frac{\pa B^3}{\pa q^3}=0.
\ee
As for $\{  p_i,q^j\}$  we  see  that  the  product
\be
\D  p_i\we\D
q^j\we(\vec{B}\D \vec{p})\we\D q^1\we\D q^2\we\D q^3=0,
\ee
therefore it remains unchanged, $\{ p_i,q^j\}=\delta^j_i$; it does
not depend on the magnetic field.
\par
As for quaternary bracket
\be
\{ f_1,f_2,f_3,f_4\}\zW=\D f_1\we\D f_2\we\D f_3\we\D f_4\we\zw
\ee
we have
\be
X_{p_1,p_2,p_3}=B^i\frac{\pa}{\pa q^i},
\ee
with the standard symplectic structure it would  be  zero.  Notice
that  $X_{p_1,p_2,p_3}$  is  not  an  inner  derivation.  This  is
a peculiar aspect  of the brackets  associated  with  intermediate
powers  of  $\zw$,  from  $\zw^2$  to  $\zw^{n-1}$,   i.e.   their
`hamiltonian vector fields' do  not preserve the bracket.
\par
As a  further  example  of  brackets  associated  with  powers  of
$\omega$ we consider an action of a Lie group $G$ on $T^*{\R}^3$.
If $G$ is any simple 3-dimensional Lie group acting on $\R^3$ with
the corresponding canonical action on $T^*\R^3$ and the associated
momentum map $J:T^*\R^3\ra\g^*$,  we  find  that  $X_{J_1,J_2,J_3}$
corresponds to  the  vector  field  associated  with  the  Casimir
function on $\g^*$ and, moreover, is an inner  derivation  of  the
quaternary bracket. In general we find
\be
X_{f_1,f_2,f_3}=\{ f_1,f_2\}X_{f_3}+\{ f_2,f_3\}X_{f_1}+
\{ f_3,f_1\}X_{f_2}
\ee
which explains why vector fields associated with  three  functions
are not inner  derivations.  Also  for  3-dimensional  simple  Lie
algebras  we   have   $C=J_1^2\pm   J_2^2\pm   J_3^2$   with   $\{
J_i,J_k\}=\pm\ze_{ikj}J_j.$

\medskip\noindent
\section{Dirac brackets}
Let us assume that on a symplectic manifold $(M,\zw)$ we have functions
$\zq_1,\dots,\zq_{2k}$ such that
\bea
&\D\zq_1\we\D\zq_2\we\cdots\we\D\zq_{2k}\ne 0,\\
&{\rm det}\left(\{\zq_j,\zq_k\}\right)\ne 0.
\eea
Dirac introduced a new Poisson bracket $\{\ ,\ \}_D$ by putting
\be
\{ f,g\}_D=\{ f,g\}-\{ f,\zq_i\}c_{ij}\{ \zq_j,g\},
\ee
where $(c_{ij})$ is the inverse of the matrix $(\{\zq_i,\zq_j\})$.
It is easy to see that $\zq_i$ are Casimir functions with respect to
this new bracket. We have the following.
\begin{Theo}
The Dirac bracket $\{\ ,\ \}_D$ is defined by the equation
\bea\nn
&(\D f\we\D g)\we\D\zq_1\we\D\zq_2\we\cdots\we\D\zq_{2k}\we\zw^{n-k-1}=\\
&\{ f,g\}_D\D\zq_1\we\D\zq_2\we\cdots\we\D\zq_{2k}\we\zw^{n-k}.
\eea
\end{Theo}
\begin{pf}
Denote by $X_f$ (resp. $X^D_f$) the hamiltonian vector field
with the hamiltonian function $f$ with respect to the original Poisson
bracket (resp. Dirac bracket). Since $\zq_i$'s are Casimir functions
with respect to the Dirac bracket we have $X^D_f(\zq_i)=0$ for $i=1,
\dots,2k$ and all possible $f$. In particular,
\be
i_{X^D_f}(\D\zq_1\we\cdots\we\D\zq_{2k})=0.
\ee
Moreover,
\bea\nn
i_{X^D_f}\zw=&i_{X_f}\zw-\{ f,\zq_i\} c_{ij}i_{X_{\zq_j}}\zw=\\
                &-\D f+\{ f,\zq_i\} c_{ij}\D\zq_j,
\eea
so that
\be
\Li_{X^D_f}\zw=\D(i_{X^D_f}\zw)=\D(\{ f,\zq_i\} c_{ij})\we\D\zq_j.
\ee
Hence,
\bea
&\Li_{X^D_f}(\D\zq_1\we\cdots\we\D\zq_{2k}\we\zw^{n-k})=\\
&\Li_{X^D_f}(\D\zq_1\we\cdots\we\D\zq_{2k})\we\zw^{n-k}+\nn\\
&(n-k)\D\zq_1\we\cdots\we\D\zq_{2k}\we\Li_{X^D_f}\zw\we\zw^{n-k-1}=0.\nn
\eea
Now, we can write
\beas
&\{ f,g\}_D\D\zq_1\we\cdots\we\D\zq_{2k}\we\zw^{n-k}=\\
&\Li_{X^D_f}(g\D\zq_1\we\cdots\we\D\zq_{2k}\we\zw^{n-k})=\\
&\D(i_{X^D_f}(g\D\zq_1\we\cdots\we\D\zq_{2k}\we\zw^{n-k}))=\\
&(n-k)\D(g\D\zq_1\we\cdots\we\D\zq_{2k}\we(i{X^D_f}\zw)\we\zw^{n-k-1})=\\
&-(n-k)\D(g\D f\we\D\zq_1\we\cdots\we\D\zq_{2k}\we\zw^{n-k-1})=\\
&(n-k)\D f\we\D g\we\D\zq_1\we\cdots\we\D\zq_{2k}\we\zw^{n-k-1}.
\eeas
\end{pf}
\medskip\noindent

\section{Generalization of previous brackets}

This  construction  can  be  generalized  by  replacing   exterior
products of  $\zw$  with  multivector  fields  of  even  order  on
manifolds of even or odd dimensions. Now we have  to  require  the
vanishing of the Schouten brackets because this  will  not  follow
automatically. Therefore we can apply  our  procedure  to  general
manifolds and general multivectors (see \cite{ILMD}.
\par
At the moment it is not easy  to  exhibit  applications  of  these
brackets to interesting physical systems. Definitely, we could use
them to select dynamical systems (vector fields)  on  the  carrier
space $M$ either by requiring the fields to be derivations of  the
brackets  or  by  associating  vector  fields  with  $k$-ples   of
functions. This would allow for a classification that goes  beyond
Hamiltonian or non-Hamiltonian dynamics.
\par
We would like to comment also that these  brackets  arising  on  a
symplectic manifold are all `natural' in the  given  context.  The
generalization to arbitrary manifolds and  arbitrary  multivectors
will lose many of the properties that we encounter on a symplectic
manifold, nevertheless they are worth investigating if  one  keeps
in mind possible applications for dynamical systems.

%%%%%%%%%%%%%%%%%%%%%%%%%%%%%%%%%%%%%%%%%%%%%%%%%%%%%%%%%%%%%%%%

\end{document}